\documentclass[10pt]{amsart}

\usepackage{epsfig,graphicx,times}
\usepackage{amssymb,amsmath}
\usepackage{enumerate}
\usepackage[T1]{fontenc}
\usepackage{ae} 

\newtheorem{theorem}{Theorem}
\newtheorem{lemma}{Lemma}

\newtheorem{remark}{Remark}
\newtheorem{definition}{Definition}

\newtheorem{conjecture}{Conjecture}
\newtheorem{corollary}{Corollary}

\def\Ps{\mathcal{P}}

\newcommand{\C}{\mathcal{C}}

\newcommand{\Z}{\mathbb{Z}}

\newcommand{\F}{\mathbb{F}}
\newcommand{\N}{\mathbb{N}}
\newcommand{\Paley}{\operatorname{Paley}}

\newcommand{\AGammaL}{\operatorname{A\Gamma L}}
\providecommand{\abs}[1]{\lvert#1\rvert}

\DeclareMathOperator{\Aut}{Aut}
\DeclareMathOperator{\AffAut}{AffAut}

\begin{document}

\title{Maximal integral point sets in affine planes over finite fields}

\author{Michael Kiermaier}
\author{Sascha Kurz}

\address{Department of Mathematics, University Bayreuth, D-95440 Bayreuth, Germany}
\email{\{michael.kiermaier,sascha.kurz\}@uni-bayreuth.de}


\maketitle

\begin{abstract}
  Motivated by integral point sets in the Euclidean plane, we consider integral point sets in affine planes
  over finite fields.
  An integral point set is a set of points in the affine plane $\F_q^2$ over a finite field $\F_q$, where the
  formally defined squared Euclidean distance of every pair of points is a square in $\F_q$.
  It turns out that integral point sets over $\F_q$ can also be characterized as affine point sets determining
  certain prescribed directions, which gives a relation to the work of Blokhuis.
  Furthermore, in one important sub-case integral point sets can be restated as cliques in Paley graphs of square order.

  In this article we give new results on the automorphisms of integral point sets and classify
  maximal integral point sets over $\F_q$ for $q\leq 47$.
  Furthermore, we give two series of maximal integral point sets and prove their maximality.
\end{abstract}

\section{Introduction}
\label{sect:intr}
The study of geometrical objects with integral edge lengths has been attractive for mathematicians throughout the ages. The first result may be obtained be the Pythagoreans considering boxes with integral side and diagonal length. A slight generalization of this problem remains unsolved. Is there a perfect box? This is a rectangular parallelepiped with all edges, face diagonals and space diagonals of integer lengths \cite{UPIN,unsolved_1}. In a more general context one is interested in the study of integral point sets, see \cite{integral_distances_in_point_sets,phd_kurz,paper_alfred} for an overview. As originally introduced integral point sets are sets of $n$ points in the $m$-dimensional Euclidean space $\mathbb{E}^m$ with pairwise integral distances. Here the majority of results are for dimension $m=2$, see for example \cite{integral_distances_in_point_sets,minimum_diameter,kreisel,phd_kurz,paper_alfred,note_on_integral_distances}.
Although integral point sets were studied for a long time our knowledge is still very limited.

Stancho Dimiev \cite{Dimiev-Setting} came up with the idea of studying integral point sets over finite fields in the hope that the situation in the finite case is easier and that some structure of the problem may be preserved.
So for a finite field $\F_q$ we consider point sets $\Ps\subseteq\F_q^2$.

In \cite{algo} the Bulgarian group around Dimiev considered integral point sets over $\F_p^2$ for primes $p\equiv 3\pmod 4$ which are maximal with respect to inclusion.
They classified the maximal integral point sets up to isomorphism for $p\in\{7,11\}$ and conjectured that the maximal integral point sets have either cardinality $\frac{p+3}{2}$ or $p$.
In the latter case all $p$ points are on a line.
Theorem~\ref{thm_kurz} clarifies the situation for cardinality $p$.
In this article we disprove their conjecture about the spectrum of possible cardinalities of maximal integral point sets and classify them for $q\le 47$.

It will turn out that in the case $q\equiv 3\pmod 4$, maximal integral point sets correspond to maximal cliques in Paley graphs of order $q^2$.
For this relevant case, in \cite{bakeretal} a construction of maximal cliques of size $\frac{q+3}{2}$ is given, which are believed to be maximal cliques of the second largest size.
They also remarked that there are other types of maximal cliques of size $\frac{q+3}{2}$.
We give a corrected proof for the maximality of these cliques and generalize the result to maximal integral point sets for $q\equiv 1\pmod 4$.
Furthermore, we give a second series of large maximal integral point sets and prove maximality. For $q\equiv 3\pmod 4$ the size is the same as in the first construction, and we believe that for $q$ sufficiently large, there are no other isomorphism classes of maximal integral point sets of this size.

\section{Definitions and basic properties}
In this article, $q=p^r$ will always be a power of a prime $p$, and unless otherwise noted, $p$ will be assumed to be odd.
$\F_q$ denotes the finite field of order $q$.
The set of squares in $\F_q$ including the zero element will be denoted by $\square_q$. 
When we talk about points and lines in $\F_q^2$, we refer to the affine plane over $\F_q$.

\begin{definition}
We define the \emph{norm} $N:\F_q^2\rightarrow\F_q$ of a point $x=(x_1,x_2)$ by $N(x):=x_1^2 + x_2^2$.
Two points $\{x,y\}\subseteq\F_q^2$ are said to be at \emph{integral distance}, if $N(x-y)\in\square_q$.
A point set $\Ps$ is called \emph{integral}, if all pairs of points are at integral distance.
An integral point set $\Ps$ is called \emph{maximal integral point set}, if there is no integral point set in $\F_q^2$ containing $\Ps$ properly.
\end{definition}

For two points $x$ and $y$, $N(x-y)$ can be interpreted as the squared Euclidean distance between $x$ and $y$.

For even numbers $q$ the situation is trivial:
we have $\square_q=\F_q$ which means that every point set $\Ps\subseteq\F_q^2$ is integral, and the only maximal integral point set is $\F_q^2$.
For that reason, we require $q$ to be odd and have $\abs{\square_q} = \frac{q+1}{2}$.

The above definition can be extended directly to finite rings instead of finite fields. For residue class rings $\mathbb{Z}_n$ the first results were obtained in \cite{Dimiev-Setting,axel_1}.

\begin{definition}
For a line $L = P + \F_q\cdot Q$ with $\{P,Q\}\subseteq\F_q^2$, $Q\neq (0,0)$, we call $D = \F_q\cdot Q$ the direction of $L$.
If $N(Q)\in\square_q$, then $D$ is called an \emph{integral direction}.
If $N(Q) = 0$, then $D$ is called a \emph{vanishing direction}.
\end{definition}

This is well-defined, because for all $\lambda\in\F_q^*$ we have $N(Q)\in\square_q$ if and only if $N(\lambda Q)\in\square_q$.

We see that two distinct points $x$ and $y$ are at integral distance if and only if $\F_q\cdot (x-y)$ is a integral direction.
So the set of integral directions completely determines if two arbitrary points $x$ and $y$ are integral distance.
Every line in $\F_q^2$ of integral direction is an integral point set, and a line not of integral direction does not contain two distinct points which are at integral distance.

For a point set $\Ps\subseteq\F_q^2$ a direction $D$ is called \emph{direction determined by $\Ps$}, if $D = \F_q\cdot(P-Q)$ with $\{P,Q\}\subseteq\Ps$.

The directions can also be identified with its slopes, which are elements of $\F_q\cup\{\infty\}$:
A direction $\F_q(x,y)\neq \F_q (0,1)$ has the slope $y/x$, and the direction $\F_q (0,1)$ corresponds to the slope $\infty$.

For $q$ odd we have $-1\in\square_q$ if and only if $q\equiv 1\pmod 4$. In this case $\omega$ will denote a fixed element with $\omega^2 = -1$.
In some cases it will be convenient to identify the affine plane $\F_q^2$ with the ring $\F_q[i]$ where $i$ is a root of the polynomial $X^2+1\in\F_q[X]$. With this identification, the map $N:(\F_q[i],\cdot)\rightarrow(\F_q,\cdot)$ is a monoid homomorphism.
An element $z=x+yi\in\F_q[i]$ with $\{x,y\}\subseteq\F_q$ is a zero-divisor if and only if $N(z) = 0$. 
In the case $q\equiv 3 \pmod 4$ we have $-1\notin\square_q$, so $X^2+1$ is irreducible and $\F_q[i]\cong \F_{q^2}$. We see that there are no vanishing directions in this case.
For $q\equiv 1\pmod 4$, $\F_q[i]$ is a finite ring with two nontrivial ideals, namely $I_1 := \F_q(\omega + i)$ and $I_2 := \F_q(\omega-i)$. These two ideals are exactly the vanishing directions and consist of the zero-divisors of $\F_q[i]$.

For $q$ odd there exists a point $\gamma\in\F[i]$ with $N(\gamma)\notin\square_q$.
The map $\nu:\F_q^2\rightarrow\F_q^2, x\mapsto \gamma x$ is bijective and maps vanishing directions to vanishing directions and non-vanishing integral directions to non-integral directions.
Together with the total number $q+1$ of the directions and the number of the vanishing directions we get the number of integral directions as $\frac{q+1}{2}$ for $q\equiv 3\pmod 4$ and $\frac{q+3}{2}$ for $q\equiv 1\pmod 4$.

Like for the complex numbers, for $z=x+yi\in\F_q[i]$ with $\{x,y\}\subseteq\F_q$ we will use the notation $\bar{z} := x-yi$.
Using this notation, we get the identity $N(z) = z\bar{z}$.

The following lemmas are straightforward and given without proof.

\begin{lemma}
  \label{lemma_pythagorean_triples}
  Let $c\in\F_q$ and $P_c=\{(a,b)\in\F_q^2\mid a^2+b^2=c^2\}$. For $c\neq 0$ we have\\[-3mm]
  \begin{eqnarray*}
    P_0 & = & \begin{cases}\{(t,\pm t\omega,0)\mid t\in\F_q\} & \mbox{if }q\equiv 1\pmod 4\\\{(0,0,0)\} &
    \mbox{if }q\equiv 3\pmod 4\end{cases} \\
    P_c & = & \{(\pm c,0,c)\}\cup\left\{\left(\frac{t^2-1}{t^2+1}\cdot c,\frac{2t}{t^2+1}\cdot c,c\right)\mid
    t\in\F_q^*,t^2\neq 1\right\}\\
    \lvert P_0\rvert & = & \begin{cases}2q-1 & \mbox{if }q\equiv 1\pmod 4 \\ 1 & \mbox{if }q\equiv 3
    \pmod 4\end{cases}\\
    \lvert P_c\rvert & = & \begin{cases}q-1 & \mbox{if }q\equiv 1\pmod 4\\q+1 & \mbox{if } q
    \equiv 3\pmod 4\end{cases}
  \end{eqnarray*}
\end{lemma}
%
%
%

\begin{definition}
The point set $\C_q := N^{-1}(1)\subseteq\F_q^2$ is called \emph{unit circle in $\F_q^2$}.
\end{definition}

By $\zeta$ we will denote a fixed multiplicative generator of the unit circle $\C_q$. In the following two lemmas we 
state that $\C_q$ shares some important properties with the well-known unit circle in the Euclidean
plane $\mathbb{R}^2$.

\begin{lemma}
  \label{lemma_cyclic}
  The unit circle $\C_q$, considered as a subset of $\F_q[i]$, is a cyclic subgroup of $\F_q[i]^*$. The order
of $\C_q$ is given by
\[
  \abs{\C_q} = \begin{cases} q-1 & \mbox{if }q\equiv 1\pmod 4 \\ q+1 & \mbox{if }q\equiv 3\pmod 4\end{cases}
\]
\end{lemma}


\begin{lemma}
\label{lemma:unit_circle_directions}
Let $\xi\in\C_q$ and $\mathfrak{L}$ be the set of the lines through $\xi$ which intersect $\C_q$ in another point.
Then $\mathfrak{L}$ contains exactly one line of each non-vanishing integral direction different from $\F_q i\xi$.
In particular, no $3$ points of $\C_q$ are collinear, and the directions determined by $\C_q$ are exactly the non-vanishing directions.
\end{lemma}


For $q\equiv 1\pmod 4$ there is another basis of $\F_q^2$ that yields a simple representation of the norm and the vanishing lines:
Let
\[a := \frac{1}{2}(1,\omega) \quad\mbox{and}\quad b := \frac{1}{2}(1,-\omega).\]
A point $P\in\F_q^2$ is said to be given in \emph{hyperbolic coordinates} $(\alpha,\beta)\in\F_q^2$, if $P = \alpha a + \beta b$.
Then $N(P) = \alpha\beta$.
In hyperbolic coordinates the unit circle $N^{-1}(1)$ is given by the hyperbola $\{(\alpha,\beta)\in\F_q^2 : \alpha\beta = 1\}$ and the vanishing directions are $(\F_q,0)$ and $(0,\F_q)$.

\section{The graph of integral distances}
It turns out that it is useful to consider integral points sets as cliques of certain graphs.

\begin{definition}
  For a fixed prime power $q=p^r$ we define the graph $\mathfrak{G}$ with vertex set
  $\F_q^2$, where two vertices $v$ and $w$ are adjacent if $N(v-w)\in\square_q$. So two
  different vertices are connected by an edge exactly if they are at integral
  distance. The graph $\mathfrak{G}$ will be called the \emph{graph of integral
  distances}.
\end{definition}

This graph is strongly regular with parameters
\[
(v,k,\lambda,\mu)=\left(q^2,\frac{q^2-1}{2},\frac{q^2-1}{4}-1,\frac{q^2-1}{4}\right)
\]
for $q\equiv 3\pmod 4$ and
\[
(v,k,\lambda,\mu)=\left(q^2,\frac{(q-1)(q+3)}{2},\frac{(q+1)(q+3)}{4}-3,\frac{(q+1)(q+3)}{4}\right)
\]
for $q\equiv 1\pmod 4$,

Furthermore, we recall that for a prime $n\equiv 1\pmod 4$ the Paley-graph $\Paley(n)$ is
defined as the graph with vertex set $\F_n$ where two vertices $v$ and $w$ are
adjacent if $v-w\in\square_n\backslash\{0\}$. For $q\equiv 3\pmod 4$ it holds: $\mathfrak{G}\cong\Paley(q^2)$.

%
%

In \cite{Blokhuis-1984} Aart Blokhuis determined the structure of cliques of maximal
size in Paley graphs of square order:
For $q\equiv 3\pmod 4$ a clique of maximal size in $\mathfrak{G}$ is
an affine line in $\F_q^2$. This implies that the size of a maximal integral point
set in $\F_q^2$ is $q$, and -- anticipating the definitions of the next section -- these
point sets are unique up to isomorphism.
Maximal cliques of size $\frac{q+1}{2}$ and $\frac{q-1}{2}$ in Paley graphs of
square order can be found in \cite{bakeretal}.

\section{Integral Automorphisms}
\label{sec_automorphism_group}
As for all combinatorial structures, the question of the automorphisms of integral point sets arises.
In our case, there are two natural ways to define an automorphism:

\begin{definition}
Let $\phi:\F_q^2\rightarrow\F_q^2$ be a bijection.
If $\phi$ preserves integral distances, $\phi$ is called \emph{integral automorphism}.
If additionally $\phi$ is an automorphism of the affine plane $\F_q^2$, $\phi$ is called \emph{affine integral automorphism}.
The group of the integral automorphisms will be denoted by $\Aut(\F_q^2)$, and the group of the affine integral automorphisms will be denoted by $\AffAut(\F_q^2)$.
\end{definition}

$\Aut(\F_q^2)$ is exactly the automorphism group of the graph of integral distances $\mathfrak{G}$. It is clear that $\AffAut(\F_q^2)$ is a subgroup of $\Aut(\F_q^2)$. In the trivial case $q$ even we have $\Aut(\F_q^2) = S_{\F_q^2}$ and $\AffAut(\F_q^2) = \AGammaL(2,q)$. In the Paley case $q\equiv 3\pmod 4$ the automorphisms $\Aut(\F_q^2)$ are known, see \cite{MR0117223,2006math_5252K,autpal,MR2025019}. According to the following theorem, in many cases the two notions of automorphisms coincide:

\begin{theorem}
\label{thm:aut=affaut}
If either $q\equiv 3 \pmod 4$ or $q\neq 5$, $q$ prime, then $\Aut(\F_q^2) = \AffAut(\F_q^2)$.
\end{theorem}

\begin{proof}
For the proof we anticipate the classification result of Theorem~\ref{thm_kurz}.
Let $\phi:\F_q^2\rightarrow\F_q^2$ be an integral automorphism.
\begin{enumerate}
\item Let $q\equiv 3\pmod 4$.
The lines in integral direction are the only maximal integral point sets of order $q$.
So an integral automorphism maps integral lines to integral lines.
Let $L$ be a not-integral line.
Using the map $\nu$ defined in section~\ref{sect:intr}, $\nu(L)$ is a integral line.
The integral automorphism $\nu\phi\nu^{-1}$ maps the integral line $\nu(L)$ to the integral line $\nu(\phi(L))$.
So $\phi(L)$ is a non-integral line, and we get that $\phi$ maps lines to lines.

\item Let $q$ be a prime with $q\equiv 1\pmod 4$ and $q\neq 5$.
Now all maximal integral point sets in $\F_q^2$ are given by the lines in integral direction and the \emph{squared crosses}\footnote{Remember $I_1 = \F_q(\omega + i)$ and $I_2 = \F_q(\omega - i)$. Furthermore, for $j\in\{1,2\}$ the expression $I_j^2$ should be read as $I_j\cdot I_j$ (and not as $I_j\times I_j$).} $x + (I_1^2\cup I_2^2)$ and $x + \tau(I_1^2\cup I_2^2)$ with $x\in\F_q^2$ and $\tau\in\F_q\setminus\square_q$.
In a squared cross no more than $2$ points are collinear in non-vanishing direction.
Each \emph{cross} $C(x) = x + (I_1 \cup I_2)$ with $x\in\F_q^2$ contains four maximal integral point sets, namely the two vanishing lines $x+I_1$ and $x+I_2$ and the two squared crosses $x+(I_1^2\cup I_2^2)$ and $x+\tau(I_1^2\cup I_2^2)$.
One checks that the crosses are the only point sets of cardinality $2q-1$ which contain exactly four maximal integral point sets, so $\phi$ maps crosses to crosses.
Since every line in vanishing direction and every squared cross is contained in a cross, we get that $\phi$ maps non-vanishing integral lines to non-vanishing integral lines.
Like in the first case, it follows that $\phi$ maps non-integral lines to non-integral lines.
Now let $D \subseteq C(x)$ be a squared cross.
We pick two points $\{P,Q\}\subseteq D\setminus\{x\}$ with $P \in x+I_1$ and $Q\in x+I_2$.
Then $P\neq Q$ and the line $L$ through $P$ and $Q$ is a line in non-vanishing direction which intersects $D$ in exactly $2$ points.
So $\phi(L)$ is a line in non-vanishing direction which intersects the maximal integral point set $\phi(D)$ of size $q$ in exactly $2$ points.
Hence $\phi(D)$ is a squared cross and $\phi$ maps squared crosses to squared crosses.
It follows that $\phi$ maps vanishing lines to vanishing lines.
\end{enumerate}
\end{proof}

\begin{theorem}(Kurz, 2007 \cite{sascha_1})\\
  \label{thm_automorphism_old}
  Let $q\notin\{5,9\}$. Then $\AffAut(\F_q^2)$, written as mappings $\F_q[i]\rightarrow \F_q[i]$, is generated by
  \begin{enumerate}[(a)]
    \item The translations $x\mapsto x+v$ for all $v\in\F_q[i]$,
    \item The reflection $x\mapsto\bar{x}$,
    \item The spiral collineations $x\mapsto\alpha x$ for all $\alpha$ with $N(\alpha)\in\square_q\setminus\{0\}$, 
    \item $a+bi\mapsto \sigma(a) + \sigma(b)i$ for all field automorphisms $\sigma$ of $\F_q$.
  \end{enumerate}
\end{theorem}

With the help of Lemma~\ref{lemma_pythagorean_triples} we easily deduce the order of $\AffAut(\F_q^2)$ as
\[
\lvert \AffAut(\F_q^2)\rvert= \begin{cases}q^2(q-1)^2r & \mbox{if }q\equiv 1\pmod 4, q\notin\{5,9\},\\q^2(q-1)(q+1)r & \mbox{if } q\equiv 3\pmod 4.\end{cases}
\]

From Theorem~\ref{thm_automorphism_old} we can deduce:
\begin{corollary}
  \label{cor_predescribe}
  $\AffAut(\F_q^2)$ acts transitively on
  \begin{enumerate}[(a)]
  \item the pairs of points at integral non-zero distance,
  \item the pairs of points at zero distance, and
  \item the pairs of points at non-integral distance.
  \end{enumerate}
\end{corollary}

We computationally investigated the sporadic cases $q=5$ and $q=9$:
\begin{theorem}
Let $G$ be the group generated by the elements in Theorem~\ref{thm_automorphism_old}.
\begin{enumerate}[(a)]
\item $[\AffAut(\F_5^2):G]=2$, $\abs{\AffAut(\F_5^2)} = 800$ and $\abs{\Aut(\F_5^2)} = 28800$.
\item $[\AffAut(\F_9^2):G]=3$, $\abs{\AffAut(\F_9^2)} = 31104$ and $\abs{\Aut(\F_9^2)}=186624$.
\end{enumerate}
\end{theorem}

There remain the cases with $q\equiv 1\pmod 4$ and $r\geq 2$.
We want to show that for these cases $\Aut(\F_q^2)$ and $\AffAut(\F_q^2)$ do not coincide.
Using hyperbolic coordinates, we define the map $\psi: \F_q^2\rightarrow\F_q^2, (\alpha,\beta)\mapsto(\alpha,\beta^p)$.

\begin{theorem}
Let $q\equiv 1\pmod 4$ and $r\geq 2$.
The above defined map $\psi$ is an integral automorphism of $\F_q^2$, but not an affine integral automorphism.
\end{theorem}

\begin{proof}
The proof is done using hyperbolic coordinates.

We have $N(\psi(\alpha,\beta)) = \alpha\beta^p = \beta^{p-1}N(\alpha,\beta)$.
Because of $p$ odd, $\beta^{p-1}\in\square_q$ and thus $N(\psi(\alpha,\beta))\in\square_q$ if and only if $N(\psi(\alpha,\beta))\in\square_q$.
Furthermore, obviously $\psi$ is a bijection, so $\psi\in\Aut(\F_q^2)$.

Let $P = \lambda(1,1)$ with $\lambda\in\F_q$ be an arbitrary point on the line $L = \F_q (1,1)$.
Then $\psi(P) = (\lambda,\lambda^p)$, which is on $L$ if and only if $\lambda\in\F_p$.
So $\psi(L)$ is not a line and therefore $\psi\notin\AffAut(\F_q^2)$.
\end{proof}

We remark that $\psi$ is a $\F_p$-linear map of order $r$.
Some computations with \texttt{nauty} \cite{nauty} and \texttt{Magma} \cite{Magma1997} showed that for $q\in\{25,49,81,121,125,169\}$ the following conjecture is true.
Unfortunately, we were not able to prove this in general.

\begin{conjecture}
Let $q\equiv 1\pmod 4$, $r\geq 2$ and $q\neq 9$.
Then
\[\Aut(\F_q^2) = \left<\AffAut(\F_q^2),\psi\right>\quad\mbox{and}\quad\abs{\Aut(\F_q^2)} = ((q-1)qr)^2.\]
\end{conjecture}

\section{Maximal integral point sets over $\F_q^2$}
\label{section:maximal}
We want to classify maximal integral point sets in $\F_q^2$ up to integral automorphisms.

\subsection{Known classification results}
\begin{theorem}(Lov\'asz and Schrijver, 1981 \cite{ls}; Blokhuis, 1984 \cite{Blokhuis-1984}; see also Kurz, 2007 \cite{sascha_1})\\
  \label{thm_kurz}
  The maximal cardinality of an integral point set over $\F_q^2$ is given by $q$.
  If $q\equiv 3\pmod 4$ then each integral point set of maximal cardinality is isomorphic to a line in
  integral direction.
  If $q\equiv 1\pmod 4$, $q$ prime and $q\neq 5$, then each integral point set of maximal cardinality is
  isomorphic to either a line in non-vanishing integral direction, or to a line in vanishing direction,
  or to a squared cross $I_1^2\cup I_2^2$.
\end{theorem}


We mention that for the above cases, $\Aut$ and $\AffAut$ are the same groups, so it is not necessary to further specify the considered automorphisms.

The key ingredient for the proof of Theorem~\ref{thm_kurz} was a theorem on point sets over $\F_q^2$ with few directions.

\begin{theorem}{(Ball, Blokhuis, Brouwer, Storme, Sz\H{o}nyi, 1999 \cite{blokhuis}; Ball 2003 \cite{1045.51004})}\\ 
  \label{thm_directions}
  Let $f:\F_q\rightarrow\F_q$,
  where $q=p^n$, $p$ prime, $f(0)=0$. Let $N=\abs{D_f}$, where $D_f$ is the set of directions
  determined by the function $f$. 
  Let $e$ (with $0\le e\le n$) be the largest integer such that each line with slope in $D_f$ meets
  the graph of $f$ in a multiple of $p^e$ points. Then we have the following:
  \begin{enumerate}
    \item[(1)] $e=0$ and $\frac{q+3}{2}\le N\le q+1$,
    \item[(2)] $p^e\ge 2$, $e\mid n$, and $\frac{q}{p^e}+1\le N\le\frac{q-1}{p^e-1}$,
    \item[(3)] $e=n$ and $N=1$.
  \end{enumerate}
  Moreover, if $p^e>2$, then $f$ is a linear map on $\F_q$ viewed as a
  vector space over $\F_{p^e}$. If $e=0$, $N=\frac{q+3}{2}$, and $n\le 2$ then $f$ is affinely equivalent to
  $f(x)=x^{\frac{q-1}{2}}$. (All possibilities for $N$ can be determined in principle.)
\end{theorem}

Looking at Theorem~\ref{thm_directions} and the proof of Theorem~\ref{thm_kurz}, it is clear that in the missing case $q\equiv 1\pmod 4$, $r\geq 2$ all the remaining maximal integral point sets of maximal cardinality are $\F_p$-vector spaces.

\subsection{Computer search}

For the classification of maximal integral point sets up to integral automorphisms we used different algorithms, because this way we could compare the results for correctness.

In all cases we built the graph of integral distances.
Now up to isomorphism all possibilities for the first few points were determined, either by Lemma~\ref{cor_predescribe} or with the help of the program $\texttt{nauty}$ \cite{nauty}.
Using these starting configurations, we searched for all maximal cliques either using the clique-search program \texttt{cliquer} \cite{cliquer}, or a fast depth-first search.
In the latter case, the search was combined with the algorithm \cite{royle} for the rejection of isomorphic copies; in fact we used an adopted version of the program in \cite{Honold-Kiermaier-2006}.
It was clear that for efficiency reasons, isomorph rejection only makes sense at lower search levels.
So for the determination of the number of isomorphism classes, we had to eliminate isomorphic copies in the end.
For that purpose we wrote a canonizer algorithm.
The canonizer algorithm implements a mapping $\kappa: 2^{\F_q^2} \rightarrow 2^{\F_q^2}$ which is constant on the orbits of $\Aut(\F_q^2)$ acting on the point sets $2^{\F_q^2}$, and for each point set $\Ps$ the image $\kappa(\Ps)$ lies within the same orbit $\Aut(\F_q^2)(\Ps)$.
For each point set $\Ps$ the image $\kappa(\Ps)$ is called \emph{canonical representative} of the orbit $\Aut(\F_q^2)(\Ps)$.

We spent a lot of work on assuring the correctness of all the algorithms: For many cases we compared the results of the cliquer-search with our depth-first search.
Then we compared the number of isomorphism classes reported by the Royle-algorithm in its pure full isomorph rejection form with the number of isomorphism classes we got by the depth-first search with separate canonization in the end.

For the hardest cases, we tried out which components are the most efficient. For the full classification, that turned out to be our depth-first search combined with isomorph rejection on the first few search levels and separate final canonization.
This way we succeeded to classify all maximal integral point sets up to $q=47$.
By $a(q)$ we denote the total number of $\Aut(\F_q^2)$-isomorphism classes of maximal integral point sets $\Ps\subseteq\F_q^2$.
If we additionally demand $\abs{\Ps} = n$, we denote that number of isomorphism classes by $a(q,n)$.
Of course, $a(q) = \sum_{n\in\N} a(q,n)$.

Table~\ref{table:classification_1} shows our classification results for $q\equiv 1\pmod 4$, and Table~\ref{table:classification_3} for $q\equiv 3\pmod 4$.
The first line of these tables contains the values of $q$ and the second line the total number $a(q)$ of $\Aut(\F_q^2)$-isomorphism classes of maximal integral point sets.
The part below splits up this information into the numbers of isomorphism classes of maximal integral point sets of a certain size $n$:
In the first column for each row the size $n$ is given, and the rest of the columns show the corresponding values $a(q,n)$.

\setlength{\arraycolsep}{2\arraycolsep}
\def\arraystretch{1} 
\begin{table}
\begin{center}
$
\begin{array}{r|rrrrrrrr}
     q& 5& 9& 13&  17&  25&   29&     37&     41 \\
\hline
\Sigma& 1& 4& 30& 107& 488& 9693& 103604& 347761 \\
\hline
     5& 1&  &   &    &    &     &       &        \\
     6&  & 2&  2&    &    &     &       &        \\
     7&  &  & 11&   8&   9&    6&      1&      1 \\
     8&  &  &  8&  57& 122&  893&    314&   1169 \\
     9&  & 2&  5&  24& 148& 4264&  17485&  61940 \\
    10&  &  &  1&  12& 108& 2864&  44952& 149839 \\
    11&  &  &   &   2&  41& 1230&  24067&  86159 \\
    12&  &  &   &    &  23&  284&  10645&  33941 \\
    13&  &  &  3&   1&  17&  116&   4835&  10854 \\
    14&  &  &   &    &   8&   22&    906&   2891 \\
    15&  &  &   &    &   4&    6&    234&    646 \\
    16&  &  &   &    &   1&    3&     89&    136 \\
    17&  &  &   &   3&   2&    2&     55&    131 \\
    18&  &  &   &    &    &     &     11&     27 \\
    19&  &  &   &    &   1&     &      2&     16 \\
    20&  &  &   &    &    &     &      3&        \\
    21&  &  &   &    &    &     &      1&      4 \\
    22&  &  &   &    &    &     &       &      3 \\
    23&  &  &   &    &    &     &       &      1 \\
    25&  &  &   &    &   4&     &      1&      1 \\
    29&  &  &   &    &    &    3&       &        \\
    37&  &  &   &    &    &     &      3&        \\
    41&  &  &   &    &    &     &       &      3
\end{array}
$
\end{center}
\caption{Numbers of isomorphism classes of maximal integral point sets for $q\equiv 1\pmod 4$.}
\label{table:classification_1}
\end{table}

\begin{table}
\begin{center}
$
\begin{array}{r|rrrrrrrrr}
     q& 3& 7& 11& 19&  23&  27&   31&     43&     47 \\
\hline
\Sigma& 1& 2&  4& 54& 294& 645& 6005& 231890& 805783 \\
\hline
     3& 1&  &   &   &    &    &     &       &        \\
     5&  & 1&   &   &    &    &     &       &        \\
     7&  & 1&  3& 25&  85&  27&   60&     15&     12 \\
     8&  &  &   &  7& 108& 411& 2004&   1748&   1097 \\
     9&  &  &   & 19&  80& 142& 2734&  54700& 125545 \\
    10&  &  &   &   &   7&  50&  933& 109127& 434029 \\
    11&  &  &  1&  4&   9&  12&  199&  54759& 210725 \\
    12&  &  &   &   &    &    &   26&   9785&  28533 \\
    13&  &  &   &   &   4&    &   46&   1490&   4904 \\
    14&  &  &   &   &    &    &     &    156&    628 \\
    15&  &  &   &   &    &   2&     &     87&    230 \\
    16&  &  &   &   &    &    &     &       &     27 \\
    17&  &  &   &   &    &    &    2&     20&     50 \\
    19&  &  &   &  1&    &    &     &       &        \\
    23&  &  &   &   &   1&    &     &      2&        \\
    25&  &  &   &   &    &    &     &       &      2 \\
    27&  &  &   &   &    &   1&     &       &        \\
    31&  &  &   &   &    &    &    1&       &        \\
    43&  &  &   &   &    &    &     &      1&        \\
    47&  &  &   &   &    &    &     &       &      1
\end{array}
$
\end{center}
\caption{Numbers of isomorphism classes of maximal integral point sets for $q\equiv 3\pmod 4$.}
\label{table:classification_3}
\end{table}
\setlength{\arraycolsep}{0.5\arraycolsep}
\def\arraystretch{1.5} 

Clearly, the spectrum of possible cardinalities of maximal integral point sets over $\F_q^2$ is more complicated as conjectured in \cite{algo}.
Looking at the results, there are some striking regularities.
Some of them follow from Theorem~\ref{thm_kurz}, others will be formulated as conjectures. 

\begin{corollary}[to Theorem~\ref{thm:aut=affaut} and Theorem~\ref{thm_kurz}]
$a(q,n) = 0$ for $n>q$.
For $q\equiv 3\pmod 4$, $a(q,q) = 1$. 
For $q\equiv 1\pmod 4$ and $q\geq 13$, $a(q,q) \geq 3$.
If additionally $q$ is prime, $a(q,q) = 3$.
\end{corollary}

\begin{conjecture}
\label{conj:n_even}
Let $q\notin\{9,13\}$. If $a(q,n) > 0$ for $n$ even, then $a(q,n+1) > 0$.
\end{conjecture}

We are interested in the minimal size $l(q)$ of a maximal integral point set in $\F_q^2$.
Our classification results show $l(q)=7$ for $11\le q\le 47$, $q\neq 13$.
Dropping the complete classification, it was possible to compute $l(q)$ for larger $q$.
We got $l(q)=8$ for $q\in\{49,53,59,61,67,73\}$, $l(71)=9$, and $l(79)\in\{8,9\}$.

\begin{lemma}
  \label{lemma_min_max}
  Let $q\ge 5$. A maximal integral point set over $\F_q^2$ contains at least $5$ points.
\end{lemma}
\begin{proof}
  Assume that $\Ps$ is an integral point set with $\abs{\Ps}\le 4$.
  If $\Ps$ contains $\abs{\Ps}-1$ collinear points, then $\Ps$ is isomorphic to a subset either of an
  integral line, or of a squared cross, or of $\Ps_L$ (see section
  \ref{section:maximal}.\ref{subsection:lineconstruction}).
  So $\Ps$ is not a maximal integral point set.
  It remains the case $\abs{\Ps} = 4$ where $\Ps$ contains no collinear triples.
  Let $\Ps=\{P_1,P_2,P_3,P_4\}$.
  We consider the three pairs of integral lines $(P_1P_2,P_3P_4)$, $(P_1P_3,P_2P_4)$ and $(P_1P_4,P_2P_3)$.
  Because $q$ is odd, in at least one of these pairs the lines intersect in a point $P_5$.
  Since no three points of $\Ps$ are collinear, it holds $P_5\notin\Ps$.
  So $\Ps$ can be extended to the integral point set $\Ps\cup\{P_5\}$.
\end{proof}

For the proofs of Theorems~\ref{thm_ref} and~\ref{theorem:line_construction} we need the following theorem, which is a consequence of the result \cite[Theorem~5.41]{lidl} on Weil sums.

\begin{theorem}(\cite[Lemma 1]{Szonyi_blocking_sets})
  \label{thm_lidl}
  Suppose $f_1,\dots,f_k,f_{k+1},\dots,f_l$ are polynomials over $\F_q$, $q$ odd. Let $N$ denote the number of
  solutions for the following requirements:
  \begin{itemize}
    \item $f_i(x)$ is a square for $i\in\{1,\dots,k\}$;
    \item $f_i(x)$ is a non-square for $i\in\{k+1,\dots,l\}$.
  \end{itemize}
  Then
  \[
    \left|N-\frac{q}{2^l}\right|\le \frac{\sqrt{q}+1}{2}\sum_{i=1}^l \operatorname{deg}\!\left(f_i\right)
  \]
  holds, unless the product of some of the $f_i$s is constant times the square of a polynomial (in this
  case the requirements can be contradicting).
\end{theorem}

\begin{theorem}
  \label{thm_ref}
  For $\varepsilon>0$ there exists a constant $q_0$ such that for all $q\ge q_0$ it holds:
  \[l(q)\ge \left( \frac{1}{2}-\varepsilon\right)\log q.\]
\end{theorem}
\begin{proof}
  Let $\varepsilon > 0$ and let $\mathcal{P} = \left\{a_i \mid i \in\{1,\ldots,k\}\right\}$ with $a_i=(x_i,y_i)\in\F_q^2$ be an integral point set of size $k < \left(\frac{1}{2}-\varepsilon\right)\log q$.
  Then $\mathcal{P}\ne\F_q^2$ and by possibly translating $\mathcal{P}$ we may assume $(0,0)\notin\mathcal{P}$.
  For $u\in\F_q$ we define the line $L_u=\{(t,ut)\mid t\in\F_q\}$.
  Our strategy is to find a suitable $u$ such that the integral point set $\mathcal{P}$ can be enlarged by a point on $L_u$.
  Then $\mathcal{P}$ is not a maximal integral point set and the proof is complete.

  We define the polynomials $f_i:=(x_i-t)^2+(y_i-ut)^2\in\F_q[t]$ which describe the squared distances between the point $a_i$ and the points $(t,ut)$ on $L_u$.
  If $u\notin\{\pm\omega\}$ the polynomials $f_i$ have degree 2.
  A polynomial $f_i$ has a multiple root if and only if $(x_i,y_i)\in L_u$, and 
  for given indices $i\neq j$ there is at most one value of $u$ such that $f_i = f_j$.
  Thus by forbidding at most $\binom{k}{2}+k+2$ values for $u$ we can ensure that $L_u\cap\mathcal{P}=\emptyset$ and the $f_i$ are pairwise distinct polynomials of degree 2 without multiple roots.

  For $q\equiv 3\pmod 4$ the $f_i$ do not have a root in $\F_q$, since $L_u\cap\mathcal{P}=\emptyset$.
  In the case $q\equiv 1\pmod 4$ the set of all points with squared distance $0$ to a point $a_i$ is given by $S_i = a_i+\F_q (1,\pm \omega)$.
  So for $N(a_i-b_i)\neq 0$ it holds $\abs{S_i\cap S_j} = 2$, and by further forbidding at most $2\binom{k}{2}$ values for $u$ we achieve that two polynomials $f_i$ and $f_j$ do not have a common root in $\F_q$ unless $N(a_i-a_j)=0$.

  Thus after all the restrictions there remain at least $q-3\binom{k}{2}-k-2$ possibilities for the choice of $u$.
  There is a $q_1\in\N$ such that for all $q \geq q_1$ this is a positive number.
  In the following we assume $q \geq q_1$ and $u$ is a fixed admissible value.

  Let $I = \{1,\ldots,k\}$ and $X := \{x\in\F_q\mid f_i(x)\neq 0\mbox{ for all }i\in I\}$.
  To finish the proof we will show that there is a $t_0\in X$ such that
  \begin{equation}
      f_i(t_0)\in\square_q\mbox{ for all }i\in I
      \label{eq_sq}
  \end{equation}

  Now assume that there is a nonempty set $J\subseteq I$ of indices such that
  \begin{equation}
      \prod_{i\in J}f_i = \sigma g^2
      \label{eq_prod}
  \end{equation}
  with $\sigma\in\F_q$ and $g\in\F_q[t]$.
  We consider the graph $G$ on the vertex set $\{1,\ldots,k\}$, where two vertices $i\neq j$ are connected by an edge if and only if $f_i$ and $f_j$ have a common root in $\F_q$.
  According to our requirements on $u$ this implies $N(a_i-a_j) = 0$, so the line through $a_i$ and $a_j$ has either direction $(1,\omega)$ or $(1,-\omega)$. 
  In the first case, we color the edge $\{i,j\}$ red, in the second case black.
  Condition~(\ref{eq_prod}) and the conditions on $u$ imply that the vertices $J$ are on a union of cycles of alternating edge color in $G$.
  Thus $\abs{J}$ is even and $\sigma = (1 + u^2)^{\abs{J}} \in\square_q\setminus \{0\}$.

  Let $i_0\in J$ and $t_0\in X$. If $f_i(t_0)\in\square_q$ for all $i\in J\setminus\{i_0\}$, then because of~(\ref{eq_prod}) and $f_i(t_0)\neq 0$ also $f_{i_0}(t_0)\in\square_q$.
  Therefore to find a $t_0\in X$ with $f_i(t_0)\in\square_q$ for all $i\in I$, it suffices to consider only the indices in $I\setminus\{i_0\}$.
  We repeat this reduction process until we obtain an index set $I_0 \subseteq I$ such that there is no nonempty $J\subseteq I_0$ fulfilling condition~(\ref{eq_prod}).
  A $t_0\in X$ is a solution of system~(\ref{eq_sq}) if and only if it is a solution of
  \begin{equation}
      f_i(t_0)\in\square_q\mbox{ for all }i\in I_0.
      \label{eq_sq2}
  \end{equation}

  Because of 
  \[
    \lim_{q\rightarrow\infty}\frac{(\sqrt{q}+3)\log q}{q^{1/2 + \varepsilon}} = 0
  \]
  there is a $q_0\in\N$, $q_0\geq q_1$ such that for all $q \geq q_0$ it holds
  \[
    k 2^k < \left(\frac{1}{2}-\varepsilon\right)q^{1/2-\varepsilon}\log q< \frac{q}{\sqrt{q}+3}
  \]
  From that follows $\frac{q}{2^{\abs{I_0}}} - \frac{\sqrt{q}+1}{2}\cdot 2\abs{I_0} > 2k$, so according to Theorem~\ref{thm_lidl} there are more than $2k$ values $t^\prime\in\F_q$ such that $f_i(t^\prime)\in\square_q$ for all $i\in I_0$.
  Because of $\abs{\F_q\setminus X} \leq 2k$, it follows that there is a solution $t_0\in X$ of system~(\ref{eq_sq2}).
  That completes the proof.
\end{proof}

\begin{remark}
  Condition~(\ref{eq_prod}) occurs for example for $\mathcal{P} =\{\pm(1,0),\pm(0,\omega)\}$.
  Here we get 
  \[
    f_1f_2f_3f_4=\Big(((1 + u^2)t^2 + 2t + 1)((1 + u^2)t^2 - 2t + 1)\Big)^2.
  \]
\end{remark}

An immediate consequence of Theorem~\ref{thm_ref} is:
\begin{corollary}
   \[
   \lim_{q\rightarrow\infty} l(q) = \infty.
   \]
\end{corollary}

\begin{conjecture}
\label{conj:q_3}
Let $q\equiv 3\pmod 4$. We define $r(q) = 2\lfloor\frac{q+1}{8}\rfloor + 5$.
\begin{enumerate}[(a)]
\item \label{conj:q_3:continuous} $a(q,n)>0$ for $n\in\{l(q),\ldots,r(q)\}$.
\item \label{conj:q_3:discrete} Let $n\geq r(q)$. If $a(q,n) > 0$, then $a(q,n+1) = 0$.
\item \label{conj:q_3:second_largest} For $q\neq 3$ there is a second largest maximal integral point set. Its cardinality is $\frac{q+3}{2}$.
\item \label{conj:q_3:2second_largest} For $q\geq 27$, $a\left(q,\frac{q+3}{2}\right) = 2$
\end{enumerate}
\end{conjecture}

\begin{remark}
Together with conjecture~\ref{conj:n_even}, part~(\ref{conj:q_3:discrete}) of the last conjecture implies $a(q,n) = 0$ for $n>r(q)$, $n$ even.

Part~(\ref{conj:q_3:second_largest}) was already conjectured in \cite{bakeretal}.
Using \texttt{cliquer}, we verified this for all $q< 200$, $q\equiv 3\pmod 4$.
\end{remark}

Now we give two series of large maximal integral point sets.
For $q\equiv 3\pmod 4$, $q\geq 27$, we believe that these are the only two types of second largest maximal integral point sets.

\subsection{A construction based on a circle}

We define $\Ps_C := \C_q^2 \cup\{0\}$.
Geometrically this is the squared unit circle together with the origin -- a special affinely regular
$\frac{q\pm 1}{2}$-gon, see \cite{regular} for a survey. By Lemma~\ref{lemma_cyclic} we have
\[
\abs{\Ps_C} = \begin{cases}\frac{q+1}{2} & \mbox{if } q\equiv 1\pmod 4,\\\frac{q+3}{2} & \mbox{if } q\equiv 3\pmod 4.\end{cases}
\]

\begin{lemma}
\label{lemma:distance_to_1}
Let $\xi\in\C_q$, $\xi = u + vi$ with $\{u,v\}\subseteq\F_q$.
\begin{enumerate}[(a)]
\item \label{lemma:distance_to_1:1} $N(1-\xi) = 2(1-u)$
\item \label{lemma:distance_to_1:2} $N(1-\xi^2) = (2v)^2$
\end{enumerate}
\end{lemma}

\begin{proof}
$\xi\in\C_q$, so $N(\xi) = u^2+v^2 = 1$ and $N(1-\xi) = (1-u)^2 + v^2 = 2(1-u)$ which shows part~(\ref{lemma:distance_to_1:1}).
For part~(\ref{lemma:distance_to_1:2}) we have $\xi^2 = (u^2 - v^2) + (2uv)i$, so part~(\ref{lemma:distance_to_1:1}) gives $N(1-\xi^2) = 2((1-u^2)+v^2) = (2v)^2$.
\end{proof}

\begin{theorem}
  \label{theorem:circle_construction}
  $\Ps_C$ is an integral point set. For $q\notin\{5,9\}$, $\Ps_C$ is a maximal integral point set.
\end{theorem}

\begin{proof}
    We identify the affine plane $\F_q^2$ with the field $\F_q[i]$.
    The map $\rho:z\mapsto\zeta^2 z$ is an integral automorphism.
    Together with Lemma~\ref{lemma:distance_to_1}~(\ref{lemma:distance_to_1:2}) it follows that $\C_q$
    is an integral point set.
    Moreover, for all $a\in\C_q^2$, $N(a-0) = 1\in\square_q$, so $\Ps_C$ is an integral point set.

    Assume that there is a point $a\in\F_q[i]\setminus\Ps_C$ such that $\Ps_C\cup\{a\}$ is an integral point set.
    Then $N(a)\in\square_q$, let $a = t^2$ with $t\in\F_q$.
    Let $A := \C_q^2 a$.
    For $t\neq 0$ it holds $A = t \C_q^2$ or $A = \zeta t\C_q^2$.
    $\Ps':=\Ps_C\cup A$ is an integral point set, because $N(\zeta^{2n} a - \zeta^{2m} a)=
    N(a) N(\zeta^{2n} - \zeta^{2m})\in\square_q$ and $N(\zeta^{2n} a - \zeta^{2m}) = N(a - \zeta^{2(m-n)})\in\square_q$.

    Now let $\xi\in\C_q\setminus\{1\}$ and write $\xi=u+vi$ with $\{u,v\}\subseteq\F_q$.
    It follows $u\neq 1$.
    With Lemma~\ref{lemma:distance_to_1}~(\ref{lemma:distance_to_1:1}) we get $N(\xi-1) = 2(1-u)\neq 0$
    and $\xi-1$ is invertible.
    Because of $a\neq 0$, it follows that $(\zeta^2)^n a = a$ is equivalent to $(\zeta^2)^n-1 = 0$,
    hence $\abs{A} = \abs{\C_q^2}$.
    For $q\equiv 3\pmod 4$ we get $\abs{\Ps_C\cup A} = q+2$, a contradiction to the maximum cardinality of an 
    integral point set, see Theorem~\ref{thm_kurz}.

    In the case $q\equiv 1\pmod 4$ we have $\abs{\Ps_C\cup A} = q$.
    Let $L_a$ be the line through $a$ and the origin.
    Lemma~\ref{lemma:unit_circle_directions} shows that no three points of $\C_q^2$ are collinear.
    If $N(a) \neq 0$, also in $A = a\C_q^2$ no three points are collinear, and so each line contains
    at most four points of $\Ps'$ different from the origin.
    If $N(a) = 0$ then $L_a$ is a line in vanishing direction, and $A$ is a subset of this line.
    Furthermore, this line does not contain another point of $\Ps'$, so $\abs{L_a\cap\Ps'} =
    \abs{A} = \frac{q-1}{2}$.
    We use the notation and the result of Theorem~\ref{thm_directions} to get a contradiction:

    First, we assume $e>0$.
    If $N(a) = 0$, then $p$ and $\abs{L_a\cap\Ps'}$ are coprime, a contradiction to $e>0$.
    So $N(a)\neq 0$. $L_a$ contains the points $0$, $a$ and $-a$.

    If $L_a$ contains no further point of $\Ps'$, we have $p=3$, $e=1$ and $A = t\zeta \C_q^2$.
    Let $\mathfrak{L}$ be the set of lines in direction $\F_q i$.
    All lines in $\mathfrak{L}$ intersect $A$ either in zero or in two points. 
    So each line which intersects $A$ contains exactly one further point of $\C_q^2\cup\{0\}$.
    But since all but two lines in $\mathfrak{L}$ intersect $\C_q^2$ either in zero or in two
    points, there are only three possibilities for this third intersection point, namely $-1$, $0$ and $1$.
    So $\abs{A}\leq 6$, $q\leq 13$ and since $p = 3$, $q\equiv 1\pmod 4$ we get $q=9$.

    If $L_a$ contains a further point $\xi$ of $\Ps'$, then $\xi\in\C_q^2$ and $L_a$ contains also
    $-\xi$, but no further point of $\Ps'$.
    So $p = 5$ and $e = 1$.
    Now each line through $1$ contains five points of $\Ps'$, so every such line must go through the
    origin, which means that there is only one such line. It follows $q=5$.

    It remains the case $e = 0$ which corresponds to case~(1) in Theorem~\ref{thm_directions}.
    Looking at the proof of Theorem~\ref{thm_kurz} one sees that for $q\notin\{5,9\}$, $\Ps'$ is
    $\AffAut(\F_q^2)$-isomorphic to $I_1^2 \cup I_2^2$.
    It follows that $\Ps'$ is the union of two collinear point sets.
    Because of $\C_q^2\subseteq\Ps'$ and no three points of $\C_q^2$ are collinear, we get
    $\abs{\C_q^2}\leq 4$, so $q\leq 9$ and since $q\equiv 1\pmod 4$ we get $q\in\{5,9\}$, a contradiction.
\end{proof}

\begin{remark}
For $q\in\{5,9\}$ the set $\Ps_C$ can be extended to an integral point set of size $q$.
\end{remark}

\subsection{A construction based on a line}
\label{subsection:lineconstruction}

We consider the line $\F_q = \F_q\cdot 1$ in $\F_q[i]$.
For each point $P\in\F_q[i]\setminus\F_q$, we define the set $S(P)\subseteq\F_q$ of all points on $\F_q$ which are at integral distance to $P$.
One of the lines in integral directions through $P$ is parallel with $\F_q$. The other lines at integral direction intersect $\F_q$, so $\abs{S(P)}$ equals the number of integral directions reduced by one.

Now we define \[
\Ps_L = S(i)\cup\{i,-i\}
\]
Geometrically, $\Ps_L$ is about the half of a line together with two points not on the line, which are arranged opposite to each other. We have
\[
\abs{\Ps_L} = \begin{cases}\frac{q+5}{2} & \mbox{if } q\equiv 1\pmod 4,\\\frac{q+3}{2} & \mbox{if } q\equiv 3\pmod 4.\end{cases}
\]

%

\begin{theorem}
  $\Ps_L$ is an integral point set. For $q\neq 9$, $\Ps_L$ is a maximal integral point set.
  \label{theorem:line_construction}
\end{theorem}
\begin{proof}
    We identify the affine plane $\F_q^2$ with the field $\F_q[i]$.
    It follows directly from the definition that $\Ps_L$ is an integral point set.
    For a point $P=x+iy\in\F_q[i]\setminus\F_q$, let $\sigma_P$ be the map $\F_q[i]\rightarrow\F_q[i]$, $z\mapsto x+yz$.
    It is easily checked that for all $P\in\F_q[i]\setminus\F_q$ we have $\sigma_P(\F_q)=\F_q$,
    $\sigma_P(i)=P$, $\sigma_P$ is an automorphism and $\sigma_P(S(i)) = S(\sigma_P(i))$.
    Now we define the set of automorphisms $G = \{\sigma_P:P\in\F_q[i]\setminus\F_q\}$.
    Clearly, $G$ is a subgroup of $\AffAut(\F_q^2)$ and operates regularly on $\F_q[i]\setminus\F_q$.

    Now assume that $\Ps_L$ is not maximal. Then there exists a point $P = (x,y)\in\F_q^2\setminus\Ps_L$
    with $S(P) = S(i)$.
    We will show that for $q\neq 9$ that leads to a contradiction.

    Let $z\in S(i)$. Then also $-z\in S(i)$ and since $S(i) = S(P)$ we get that also $2x+z\in S(i)$.
    It follows that $S(i)$ is the disjoint union of sets $[z] = z + 2\Z x$ where $z\in S(i)$.
    If $x\neq 0$ then $p$ divides the cardinality of $[z]$ for all $z\in S(i)$, so $p$ also divides
    $\abs{S(i)}$, a contradiction. Therefore $x = 0$ and $P = yi$.

    First we deal with the case $q\equiv 3\pmod 4$:
    We have $\sigma_P(z) = yz$, so $\sigma_P$ can be considered as an element of $\F_q^*$.
    Hence the order $o$ of $\sigma_P$ divides $q-1$.
    Because of $P\notin\{\pm i\}$ we have $o\geq 3$.
    We consider two cases:
    \begin{enumerate}[(1)]
    \item \label{thm:line_construction:case1} The case $o$ is prime:\\
    $\sigma_P$ has exactly one fixed point on $S(i)$, namely the origin.
    We have $\sigma_P(S(i)) = S(\sigma_P(i)) = S(P) = S(i)$, so the group action of
    $\left<\sigma_P\right>$ partitions $S(i)$ into this fixed point and orbits of length $o$.
    Hence $o$ divides $\abs{S(i)}-1 = \frac{q-3}{2}$, which gives $o = 2$, a contradiction.
    \item \label{thm:line_construction:case2} The case $o$ is composite:\\
    Because of $4\nmid q-1$, $o$ has a prime divisor $s\neq 2$.
    We define $Q := \sigma_P^{o/s}(P)$.
    Then $S(Q) = S(P) = S(i)$ and the order of $\sigma_Q = \sigma_P^{o/s}$ is $s$.
    Now we apply case~(\ref{thm:line_construction:case1}) to $Q$ and get $s = 2$, a contradiction.
    \end{enumerate}

    For the remaining case $q\equiv 1\pmod 4$, $q\neq 9$ we first consider two separate cases:
    For $q=5$ we have $\abs{\Ps_L} = 5$, so $\Ps_L$ is maximal by Theorem~\ref{thm_kurz}.
    For $q\in\{13,17\}$ we use a computer program to see that $\Ps_L$ is maximal.
    So in the following we can assume $q\geq 25$.

    Similarly as in Lemma~\ref{lemma_pythagorean_triples} we get the identity $S(i) =
    \left\{\frac{1-t^2}{2t}\mid t\in\F_q^*\right\}$.
    Because of $S(i) = S(P)$, for all $t\in\F_q^*$ the term $1+y^2\left(\frac{1-t^2}{2t}\right)^2$ is a square.
    We define the polynomial
    \[
      f=4t^2+y^2(1-t^2)^2\in\F_q[t]
    \]
    and conclude $f(t)\in\square_q$ for all $t\in\F_q$.
    So we get
    \[
    N := \abs{\{x\in\F_q\mid f(x)\in\square_q\}}=q
    \]
    For $q\geq 25$ it holds $\abs{N-\frac{q}{2}} = N-\frac{q}{2} > \deg(f)\frac{\sqrt{q}+1}{2}$, so by Theorem~\ref{thm_lidl}
    $f$ is a constant times the square of a polynomial. Thus $f$ contains a repeated factor of
    degree either $1$ or $2$.
    \begin{enumerate}[(1)]
        \item The repeated factor has degree $1$:\\
        There exists an $t_0\in\F_q$ with $f(t_0)=f'(t_0)=0$. We have
        \[
        f'(t_0)=4 t_0(2-y^2+y^2t_0^2)=0
        \]
        Because of $f(0) = y^2\neq 0$, it follows $t_0\neq 0$, so $8-4y^2+4y^2t_0^2 = 0$ which
        leads to $t_0^2 = 1-2/y^2$.
        Now from $f(t_0) = 0$ we get $y\in\{\pm 1\}$, which contradicts $P\notin\Ps_L$.
        \item The repeated factor has degree $2$:\\
        Here $f = b(t^2 + a)^2$ with $\{a,b\}\subseteq\F_q$.
        Comparing coefficients gives $b=y^2$, $a^2=1$ and $(a+1)y^2=2$.
        It follows $a = 1$ and $y\in\{\pm 1\}$, again a contradiction.
    \end{enumerate}
\end{proof}

\begin{remark}
\begin{enumerate}[(a)]
\item For $q\equiv 3\pmod 4$ the above proof was already given in \cite{bakeretal}.
Nevertheless we decided to give the proof here, because on the one hand for the case $q\equiv 1\pmod 4$ the first part of the proof including the notation was needed, and on the other hand, there is a flaw in the original proof: 
They use the argumentation of the above case~(\ref{thm:line_construction:case1}) regardless if $o$ is prime or not.
But for composite numbers $o$ the non-singleton orbits can have sizes different from $o$.
We took care of this by introducing case~(\ref{thm:line_construction:case2}) for composite $o$.
\item For $q = 9$ we have $S(i) = S(\omega i) = \{0,1,-1,\omega,-\omega\}$. So $\Ps_L$ can be extended by $\pm\omega i$ to an integral point set of maximal cardinality $9$. For a primitive element $\alpha$, the polynomial $f$ in the above proof has the form 
\[f = -(t-\alpha)(t-\alpha^3)(t-\alpha^5)(t-\alpha^7)\]
without repeated factors.
\item For $q\equiv 1\pmod 4$ one could try to do a similar construction starting with a line in vanishing direction.
But the resulting integral point set always can be extended to an integral point set of order $q$ isomorphic to $I_1^2\cup I_2^2$.
\end{enumerate}
\end{remark}

\begin{figure}[htp]
  \begin{center}
    \setlength{\unitlength}{0.20cm}
    \begin{picture}(23.2,23.2)
      \multiput(0.15,0)(0,1){24}{\line(1,0){23}}
      \multiput(0.15,0)(1,0){24}{\line(0,1){23}}
      \put(0.65,0.5){\circle*{0.6}}
      \put(0.65,22.5){\circle*{0.6}}
      \put(22.65,0.5){\circle*{0.6}}
      \put(0.65,1.5){\circle*{0.6}}
      \put(1.65,0.5){\circle*{0.6}}
      \put(4.65,19.5){\circle*{0.6}}
      \put(19.65,19.5){\circle*{0.6}}
      \put(19.65,4.5){\circle*{0.6}}
      \put(4.65,4.5){\circle*{0.6}}
      \put(0.65,14.5){\circle*{0.6}}
      \put(14.65,0.5){\circle*{0.6}}
      \put(0.65,9.5){\circle*{0.6}}
      \put(9.65,0.5){\circle*{0.6}}
    \end{picture}
    \quad
    \setlength{\unitlength}{0.20cm}
    \begin{picture}(23.2,23.2)
      \multiput(0.15,0)(0,1){24}{\line(1,0){23}}
      \multiput(0.15,0)(1,0){24}{\line(0,1){23}}
      \put(0.65,0.5){\circle*{0.6}}
      \put(11.65,15.5){\circle*{0.6}}
      \put(11.65,8.5){\circle*{0.6}}
      \put(1.65,0.5){\circle*{0.6}}
      \put(1.65,7.5){\circle*{0.6}}
      \put(21.65,0.5){\circle*{0.6}}
      \put(1.65,16.5){\circle*{0.6}}
      \put(3.65,21.5){\circle*{0.6}}
      \put(17.65,0.5){\circle*{0.6}}
      \put(3.65,2.5){\circle*{0.6}}
      \put(19.65,5.5){\circle*{0.6}}
      \put(19.65,18.5){\circle*{0.6}}
      \put(8.65,0.5){\circle*{0.6}}
    \end{picture}
  \end{center}
  \caption{The integral point sets $\Ps_1$ and $\Ps_2$.}
  \label{fig_sporadic_23}
\end{figure}

\begin{corollary}
  Let $q\equiv 3\pmod 4$ and $q\geq 11$. Then $a\left(q,\frac{q+3}{2}\right) \geq 2$.
\end{corollary}

\begin{proof}
By Theorem~\ref{theorem:circle_construction} and Theorem~\ref{theorem:line_construction} we have two constructions $\Ps_C$ and $\Ps_L$ of maximal integral point sets of size $\frac{q+3}{2}$.
For $q\geq 11$, there is a collinear subset of $\Ps_L$ of size $4$.
$\Ps_C$ does not admit a collinear subset of size $4$, so by Theorem~\ref{thm:aut=affaut}, $\Ps_C$ and $\Ps_L$ are not $\Aut(\F_q^2)$-isomorphic.
\end{proof}

For $q\in\{3,7\}$, we have $\Ps_C=\Ps_L$.
For $q\equiv 3\pmod 4$, $27\leq q\leq 47$, our classification shows that there are no other isomorphism classes of maximal integral point sets of size $\frac{q+3}{2}$, and we conjecture that this is true for all $q\geq 27$, see conjecture~\ref{conj:q_3}~(\ref{conj:q_3:2second_largest}).

But for $q\in\{11,19,23\}$ there are other such point sets, which also have a nice geometric pattern.
For $q=23$ these examples are shown in figure~\ref{fig_sporadic_23}. They are
\[
  \mathcal{P}_1=\{0\}\cup 1\cdot \C_{23}^6\cup 3\cdot \C_{23}^6\cup  9\cdot \C_{23}^6
\]
and
\[
  \mathcal{P}_2=\{0\}\cup 1\cdot \C_{23}^8\cup 2 \zeta^4\cdot\C_{23}^8\cup  6\zeta^4\cdot\C_{23}^8\cup 8\cdot\C_{23}^8
\]
For $q=19$ one of the two examples has a similar shape and is shown in figure~\ref{fig_sporadic_19_1}. It is
\[
  \mathcal{P}_3=\{0\}\cup 1\cdot\C_{19}^4\cup 3\cdot\C_{19}^4
\]
\begin{figure}[htp]
  \begin{center}
    \setlength{\unitlength}{0.20cm}
    \begin{picture}(19.2,19.2)
      \multiput(0.15,0)(0,1){20}{\line(1,0){19}}
      \multiput(0.15,0)(1,0){20}{\line(0,1){19}}
      \put(0.65,0.5){\circle*{0.6}}
      \put(7.65,3.5){\circle*{0.6}}
      \put(2.65,4.5){\circle*{0.6}}
      \put(2.65,15.5){\circle*{0.6}}
      \put(7.65,16.5){\circle*{0.6}}
      \put(1.65,0.5){\circle*{0.6}}
      \put(2.65,9.5){\circle*{0.6}}
      \put(6.65,12.5){\circle*{0.6}}
      \put(6.65,7.5){\circle*{0.6}}
      \put(2.65,10.5){\circle*{0.6}}
      \put(3.65,0.5){\circle*{0.6}}
    \end{picture}
  \end{center}
  \caption{The integral point set $\Ps_3$.}
  \label{fig_sporadic_19_1}
\end{figure}
The second sporadic example $\Ps_4$ for $q=19$ and the sporadic example $\Ps_5$ for $q=11$ have a different geometric pattern. They are subsets of $\C_q\cup\F_q$, see figure~\ref{fig_sporadic_line_and_circle}.

\begin{figure}[htp]
  \begin{center}
    \setlength{\unitlength}{0.20cm}
    \begin{picture}(20.2,20.2)
      \multiput(0.15,0)(0,1){20}{\line(1,0){19}}
      \multiput(0.15,0)(1,0){20}{\line(0,1){19}}
      \put(0.65,0.5){\circle*{0.6}}
      \put(0.65,1.5){\circle*{0.6}}
      \put(0.65,18.5){\circle*{0.6}}
      \put(4.65,0.5){\circle*{0.6}}
      \put(4.65,2.5){\circle*{0.6}}
      \put(4.65,17.5){\circle*{0.6}}
      \put(5.65,0.5){\circle*{0.6}}
      \put(14.65,0.5){\circle*{0.6}}
      \put(15.65,0.5){\circle*{0.6}}
      \put(15.65,2.5){\circle*{0.6}}
      \put(15.65,17.5){\circle*{0.6}}
    \end{picture}
    \quad 
    \setlength{\unitlength}{0.20cm}
    \begin{picture}(11.2,11.2)
      \multiput(0.15,0)(0,1){12}{\line(1,0){11}}
      \multiput(0.15,0)(1,0){12}{\line(0,1){11}}
      \put(0.65,0.5){\circle*{0.6}}
      \put(0.65,1.5){\circle*{0.6}}
      \put(0.65,10.5){\circle*{0.6}}
      \put(2.65,0.5){\circle*{0.6}}
      \put(3.65,5.5){\circle*{0.6}}
      \put(3.65,6.5){\circle*{0.6}}
      \put(6.65,0.5){\circle*{0.6}}
    \end{picture}
  \end{center}
  \caption{The integral point sets $\Ps_4$ and $\Ps_5$.}
  \label{fig_sporadic_line_and_circle}
\end{figure}

\section{Remarks on integral point sets over $\mathbf{\mathbb{E}^2}$}
It is interesting to mention that the situation for integral point sets in $\mathbb{E}^2$ is somewhat similar.
Instead for maximal integral point sets, here we ask for the minimum possible diameter $d(2,n)$ of an integral point set in the Euclidean plane $\mathbb{E}^2$ with pairwise integral distances, where the diameter is the largest occurring distance.
Without any extra condition $n$ points on a line would yield an integral point set with small diameter.
To make it more interesting one forces integral point sets in $\mathbb{E}^2$ to be two dimensional.
Here all known non-collinear examples of integral point sets with minimum diameter consist of a line with $n-1$ points and one point apart, see \cite{paper_alfred,note_on_integral_distances}.
If we forbid $3$ points to be collinear, integral point sets on circles seem to be the examples with minimum diameter.
The situation stays more or less the same if we consider integral point sets over $\mathbb{Z}^2$.
These results on the structure of integral point sets over $\mathbb{E}^2$ or $\mathbb{Z}^2$ are up to now only conjectures which are verified for the first few numbers $n$ of points. So this is one motivation to study integral point sets over $\F_q^2$ in the hope that here the situation is easier to handle.

Besides the characterization of the maximal integral point set with largest or second largest cardinality another interesting question is the characterization of those maximal integral point sets with minimum cardinality.
From our data we may conjecture that for $q\ge 11$ we have $a(q,s)=0$ for $s\le 6$.
Again we can compare this situation to the situation in $\mathbb{E}^2$.
A result due to Almering \cite{rationale_vierecke} is the following.
Given any integral triangle $\Delta$ in the plane, the set of all points $x$ with rational distances to the three corners of $\Delta$ is dense in the plane.
Later Berry generalized this results to to triangles where the squared side lengths and at least one side length are rational.
In $\mathbb{Z}^2$ the situation is a bit different.
In \cite{paper_carpet} the authors search for maximal integral triangles over $\mathbb{Z}^2$.
They exist but seem to be somewhat rare. There are only seven maximal integral triangles with diameter at most $5000$.
The smallest possible diameter is $2066$.
In a forthcoming paper \cite{heronian} one of the authors has extended this list, as a by-product, up to diameter $15000$ with in total $126$ maximal integral triangles.

\section{Conclusion}
In this article we gave new results on the automorphisms of integral point sets and we classified
maximal integral point sets over $\F_q$ for $q\leq 47$. Furthermore, we gave two series of maximal
integral point sets and have proven their maximality. There remain several conjectures to be settled.

We saw that for integral point sets over finite fields, there is a big difference between the cases $q\equiv 3\pmod 4$ and $q\equiv 1\pmod 4$.
Within the latter case, the cases $q\in\{5,9\}$ significantly differ from the others.
The former case apparently admits stronger conclusions and shorter proofs.
While some results (see Theorems~\ref{thm:aut=affaut} and~\ref{thm_automorphism_old})
do no directly generalize, others
(see Theorems~\ref{thm_ref},~\ref{theorem:circle_construction} and~\ref{theorem:line_construction}) can be carried forward but require more elaborate proofs.

The difference between the two cases arises from the fact that $x^2+y^2\in\F_q[x,y]$ is irreducible only for $q\equiv 3\pmod 4$. More generally one could consider distance functions $x^2+ky^2$. If it is irreducible the unit circle
$\left\{(x,y)\in\F_q^2\mid x^2+ky^2=1\right\}$ geometrically corresponds to
an ellipse and otherwise to a hyperbola.

\section*{Acknowledgment}
We would like to thank the anonymous referees for their very useful comments, remarks, and pointers to the literature. Especially the proof of Theorem \ref{thm_ref} is based on one of their suggestions.

\end{document}